\documentclass[12pt]{article}
\usepackage{amsfonts,amsmath,latexsym,amssymb,mathrsfs}
\usepackage{breakcites}
\usepackage[dvips]{graphicx}
\usepackage{wrapfig}
\usepackage{caption}
\usepackage{color}
\usepackage{geometry}
\usepackage{enumerate}
\usepackage[T1]{fontenc}
\usepackage[english]{babel}
\usepackage{dsfont,bm,bbm,mathrsfs}
\usepackage{stmaryrd}
\usepackage{url,color}
\usepackage{textcase}
\usepackage{multirow}
\usepackage{tikz}
\usetikzlibrary{intersections}

\allowdisplaybreaks

\newtheorem{thm}{Theorem}[section]
\newtheorem{prop}[thm]{Proposition}

\newtheorem{lma}[thm]{Lemma}
\newtheorem{cor}[thm]{Corollary}
\newtheorem{exam}[thm]{Example}

\newtheorem{rem}[thm]{Remark}
\newtheorem{con}[thm]{Conjecture}

\makeatletter
\def\be#1{\begin{equation*}#1\end{equation*}}
\def\ben#1{\begin{equation}#1\end{equation}}

\def\bea#1{\begin{align*}#1\end{align*}}
\def\bean#1{\begin{align}#1\end{align}}

\def\given{\mskip 0.5mu plus 0.25mu\vert\mskip 0.5mu plus 0.15mu}
\newcounter{@bracketlevel}
\def\@bracketfactory#1#2#3#4#5#6{
\expandafter\def\csname#1\endcsname##1{%
\addtocounter{@bracketlevel}{1}%
\global\expandafter\let\csname @middummy\alph{@bracketlevel}\endcsname\given%
\global\def\given{\mskip#5\csname#4\endcsname\vert\mskip#6}\csname#4l\endcsname#2##1\csname#4r\endcsname#3%
\global\expandafter\let\expandafter\given\csname @middummy\alph{@bracketlevel}\endcsname
\addtocounter{@bracketlevel}{-1}}%
}
\def\bracketfactory#1#2#3{%
\@bracketfactory{#1}{#2}{#3}{relax}{0.5mu plus 0.25mu}{0.5mu plus 0.15mu}
\@bracketfactory{b#1}{#2}{#3}{big}{1mu plus 0.25mu minus 0.25mu}{0.6mu plus 0.15mu minus 0.15mu}
\@bracketfactory{bb#1}{#2}{#3}{Big}{2.4mu plus 0.8mu minus 0.8mu}{1.8mu plus 0.6mu minus 0.6mu}
\@bracketfactory{bbb#1}{#2}{#3}{bigg}{3.2mu plus 1mu minus 1mu}{2.4mu plus 0.75mu minus 0.75mu}
\@bracketfactory{bbbb#1}{#2}{#3}{Bigg}{4mu plus 1mu minus 1mu}{3mu plus 0.75mu minus 0.75mu}
}
\bracketfactory{clc}{\lbrace}{\rbrace}
\bracketfactory{clr}{(}{)}
\bracketfactory{cls}{[}{]}
\bracketfactory{abs}{\lvert}{\rvert}
\bracketfactory{norm}{\Vert}{\Vert}
\bracketfactory{floor}{\lfloor}{\rfloor}
\bracketfactory{ceil}{\lceil}{\rceil}
\bracketfactory{angle}{\langle}{\rangle}

\def\ER{Erd\H{o}s-R\'enyi}

\makeatother

\newcommand{\law}{\mathscr{L}}

\DeclareMathOperator{\E}{\mathbb{E}}

\newcommand{\mean}{\E}
\newcommand{\R}{\mathbb{R}}

\newcommand{\non}{\nonumber}
\newcommand{\Z}{\mathbb{Z}}

\DeclareMathOperator{\Var}{\mathrm{Var}}
\DeclareMathOperator{\var}{\mathrm{Var}}

\usepackage{subfig}

\newcommand{\bone}{{\bf 1}}
\newcommand{\bc}{{\mathds{C}}}

\newcommand{\cf}{{\cal F}}

\newcommand{\ci}{{\cal I}}

\newcommand{\dtv}{{d_{\rm TV}}}

\def\tf{{\tilde f}}
\newcommand{\Pn}{{\rm Pn}}

\DeclareMathOperator*{\esssup}{ess\,sup}

\def\Ref#1{(\ref{#1})}
\def\a{\alpha}
\def\b{\beta}

\def\l{\lambda}

\def\sjn{\sum_{j=1}^n}

\def\dtv{d_{\mathrm{TV}}}
\def\tpj{\overline{\tau^+_j}}
\def\tmj{\overline{\tau^-_j}}
\def\tpi{\overline{\tau^+_i}}
\def\tpi1{\overline{\tau^+_{i-1}}}
\def\tmi{\overline{\tau^-_i}}

\makeatletter

\newcommand{\qed}{\nopagebreak\hspace*{\fill}
{\vrule width6pt height6ptdepth0pt}\par}

\newcommand{\red}[1]{\textcolor{black}{#1}}

\def\ignore#1{}

\makeatletter
\newcommand*\labelcounter[2]{\begingroup
  \protected@edef\@currentlabel{\csname p@#1\endcsname\csname the#1\endcsname}%
  \label{#2}\endgroup}
\newcommand*\refsetcounter[2]{\setcounter{#1}{#2}%
  \protected@edef\@currentlabel{\csname p@#1\endcsname\csname the#1\endcsname}%
  }
\makeatother

\newcounter{rtaskno}

\ignore{\newcounter{con}
}

\newcounter{con}%for constants

\numberwithin{equation}{section}

\title{\sc\bf\large\MakeUppercase{On moderate deviations in Poisson approximation}}

\author{ Qingwei Liu\footnote{School of Mathematics and Statistics,
The University of Melbourne,
VIC 3010, Australia, E-mail: qingweil@student.unimelb.edu.au. Work supported in part by China Scholarship Council.}
\ and \
 Aihua Xia\footnote{School of Mathematics and Statistics,
The University of Melbourne,
VIC 3010, Australia, E-mail: aihuaxia@unimelb.edu.au. Work supported in part by Australian Research Council Grants No DP190100613.}
}

\def\parsedate #1:20#2#3#4#5#6#7#8\empty{20#2#3-#4#5-#6#7}
\def\moddate{\expandafter\parsedate\pdffilemoddate{\jobname.tex}\empty}
\date{\moddate}

\begin{document}
\maketitle

\begin{abstract} %Poisson distribution is more suited for approximating the the counts of rare events than normal distribution.  
\red{In this paper, we first use} the distribution of the number of records to  demonstrate that the right tail probabilities of counts of rare events are generally better approximated by the right tail probabilities of Poisson distribution than {those} of normal distribution. We then show the moderate deviations in Poisson approximation generally require an adjustment and, with suitable adjustment, we establish better error estimates of the moderate deviations in Poisson approximation than those in \cite{CFS}. Our estimates contain no unspecified constants and are easy to apply. We illustrate the use of the theorems in six applications: Poisson-binomial distribution, matching problem, occupancy problem, birthday problem, random graphs and 2-runs.
The paper complements the works of \cite{CC92,BCC95,CFS}.
\end{abstract}



\section{Introduction}%\label{secIntroduction}

An exemplary moderate deviation theorem is as follows (see \cite[p.~228]{Petrov75}). Let $X_i$, $1\le i\le n$, be independent and identically distributed (i.i.d.) random variables with $\mean(X_1)=0$ and $\var(X_1)=1$. If for some $t_0>0$,
\begin{equation}
\mean e^{t_0|X_1|}\le c_0<\infty,\label{momentgeneratingcondition}
\end{equation}
then there exist positive constants $c_1$ and $c_2$ depending on $c_0$ and $t_0$ such that
\begin{equation}
\frac{\mathbb{P}\left(\frac1{\sqrt{n}}\sum_{i=1}^nX_i\ge z\right)}{1-\Phi(z)}=1+O(1)\frac{1+z^3}{\sqrt{n}},\ \ \ 0\le z\le c_1n^{1/6},\label{petrov1}
\end{equation}
where $\Phi(z)$ is the distribution function of the standard normal, $|O(1)|\le c_2$. However, since the pioneering work \cite{Chen75}, it has been shown \cite{BHJ} that, for the counts of rare events, Poisson distribution provides a better approximation. For example, the distribution of the number of records \cite{Dwass60,Renyi62} in Example~\ref{exam1} below
can be better approximated by the Poisson distribution having the same mean than by a normal distribution \cite{DP88}. \red{Moreover, a suitable} refinement of the Poisson distribution can further improve the performance of the approximation \cite{Borovkov88,BP96}.

The right tail probabilities of counts of rare events are often needed in statistical inference but these probabilities are so small that the error estimates in approximations of distributions of the counts are usually of no use because the bounds are often larger than the probabilities of interest. Hence it is of practical interest to consider their approximations via moderate deviations in Poisson approximation in a similar fashion to \Ref{petrov1}. However, there is not much progress in the general framework except the special cases in \cite{CC92,BCC95,CFS,TLX18,CV19}. This is partly due to the fact that the tail behaviour of a Poisson distribution is significantly different from that of a normal distribution and this fact is observed by \cite{Gn43} in the context of extreme value theory. In particular, \cite{Gn43} concludes that the Poisson distribution does not belong to any domain of attraction while the normal distribution belongs to the domain of attraction of the Gumbel distribution.

\begin{exam}\label{exam1} {\rm We use the distribution of the number of records to explain the difference of moderate deviations between Poisson and normal approximations.
More precisely,  let $\{\eta_i:\ 1\le i\le n\}$ be {i.i.d.} random variables with a continuous
cumulative distribution function. As the value of $\eta_1$ is always a record, for $2\le i\le n$, we say $\eta_i$ is a record if
$\eta_i>\max_{1\le j\le i-1}\eta_j$. We define the indicator random variable
$$I_i:=\bone[\eta_i>\max_{1\le j\le i-1}\eta_j],$$ 
that is, $I_i=1$ if a new record occurs at
time $i$ and $I_i=0$ otherwise. Our interest is on the distribution of $S_n:=\sum_{i=2}^nI_i$, denoted by $\law(S_n)$. \cite{Dwass60,Renyi62} state that $\mean I_i=1/i$, $\{I_i:\ 2\le i\le n\}$ are independent {so}
$$\lambda_n:=\mean S_n=\sum_{i=2}^n\frac1i;\ \ \ \sigma_n^2:=\var(S_n)=\sum_{i=2}^n\frac1i\left(1-\frac1i\right).%\label{recordmeanvar}
$$
We use $\Pn(\lambda)$ to stand for the Poisson distribution with mean $\lambda$, 
$\Pn(\lambda)(A):=\mathbb{P}(Y\in A)$ for $Y\sim\Pn(\lambda)$, and $N(\mu,\sigma^2)$ to stand for the
normal distribution with mean $\mu$ and variance $\sigma^2$.

Let $v_n:=\lambda_n+x\cdot\sigma_n$, and we consider approximations of
$\mathbb{P}(S_n\ge v_n)$ by moderate deviations based on $\Pn(\lambda_n)$  \cite{BCC95,CFS} and $N_n\sim N(\lambda_n,\sigma_n^2)$. For $x=3$, figures~\ref{figure1},~\ref{figure2}
and~\ref{figure4}
are respectively the plots of the ratios
$\mathbb{P}(S_n\ge v_n)/\Pn(\lambda_n)([v_n,\infty))$, $\mathbb{P}(S_n\ge v_n)/\mathbb{P}(N_n\ge v_n)$ and $\mathbb{P}(S_n\ge v_n)/\Pn(\sigma_n^2)([v_n,\infty))$ for the range of $n\in[3,10^5]$. 
As observed in \cite{BP96}, Poisson and normal approximations to $\law(S_n)$ are resp. with order $O((\ln n)^{-1})$ and $O((\ln n)^{-1/2})$, the numerical studies confirm that approximation by the Poisson distribution is better than that by normal distribution. In fact, it appears that the speed of convergence of 
$\mathbb{P}(S_n\ge v_n)/\mathbb{P}(N_n\ge v_n)$ to $1$ as $n\to\infty$ is too slow to be of practical use. In the context of normal approximation to the distribution of integer valued random variables, a common practice is to introduce a 0.5 correction, giving the ratios $\mathbb{P}(S_n\ge v_n)/\mathbb{P}(N_n\ge\lceil v_n\rceil-0.5)$, where $\lceil x\rceil$ is the smallest integer that is not less than $x$. Figure~\ref{figure3} is the plot of the ratios and we can see that the ratios 
are still far away from the limit of 1. Finally, the difference between Figure~\ref{figure1} 
and Figure~\ref{figure4} shows that a minor change of the mean of the approximating Poisson can change the quality of moderate deviation approximation significantly, further highlighting the difficulty of obtaining sharp bounds in theoretical studies in the area.
}

\begin{figure}[ht]
\begin{minipage}[t]{0.49\linewidth}  % <---
%the four numbers after "trim=" are for left bottom right top
   \includegraphics[trim = 10mm 75mm 5mm 80mm,height=0.26\textheight]{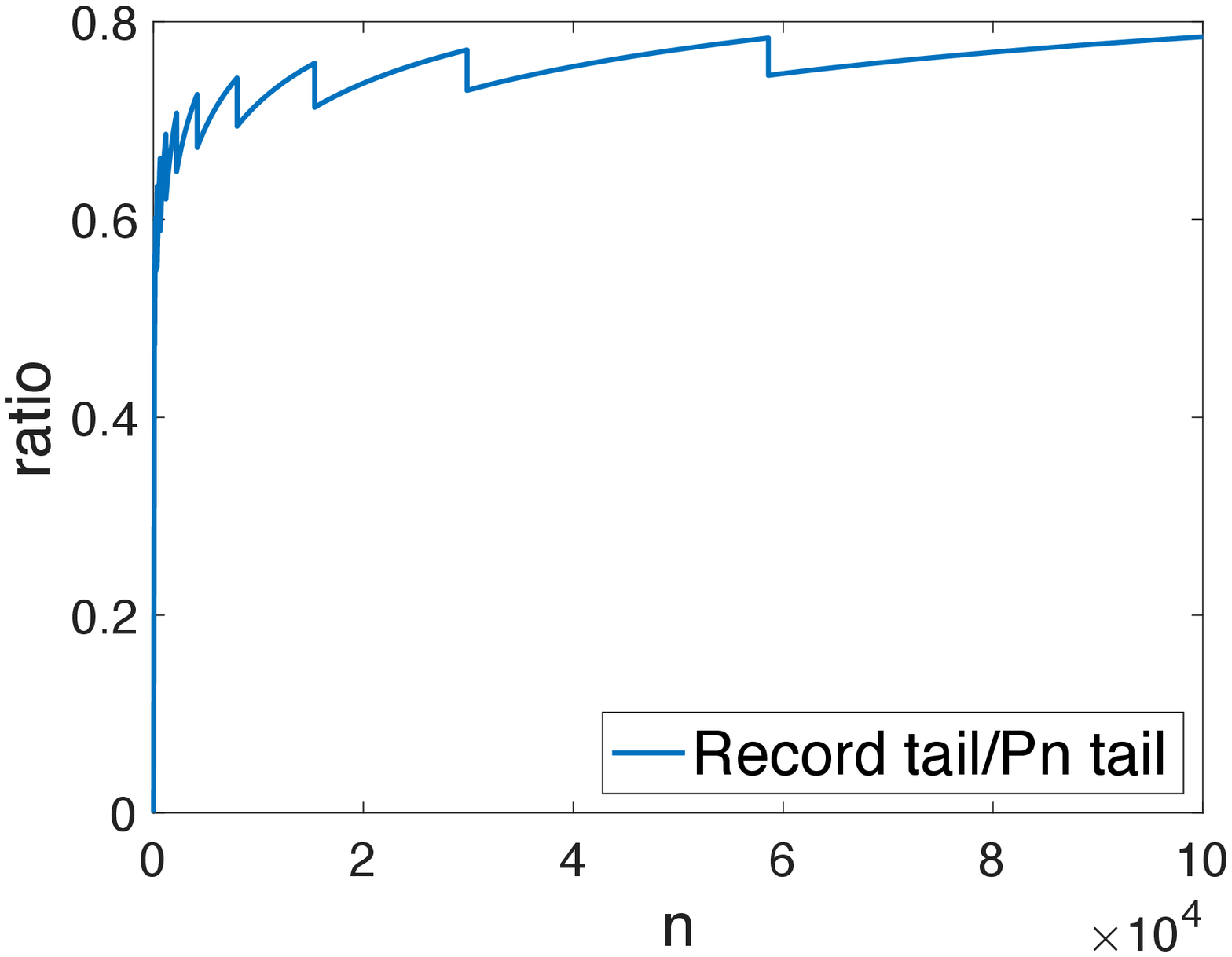}  
 \caption{Pn($\lambda_n$)} %
\label{figure1} \end{minipage}
\hfill
\begin{minipage}[t]{0.5\linewidth}  % <---
%the four numbers after "trim=" are for left bottom right top
   \includegraphics[trim = 10mm 75mm 5mm 80mm,height=0.26\textheight]{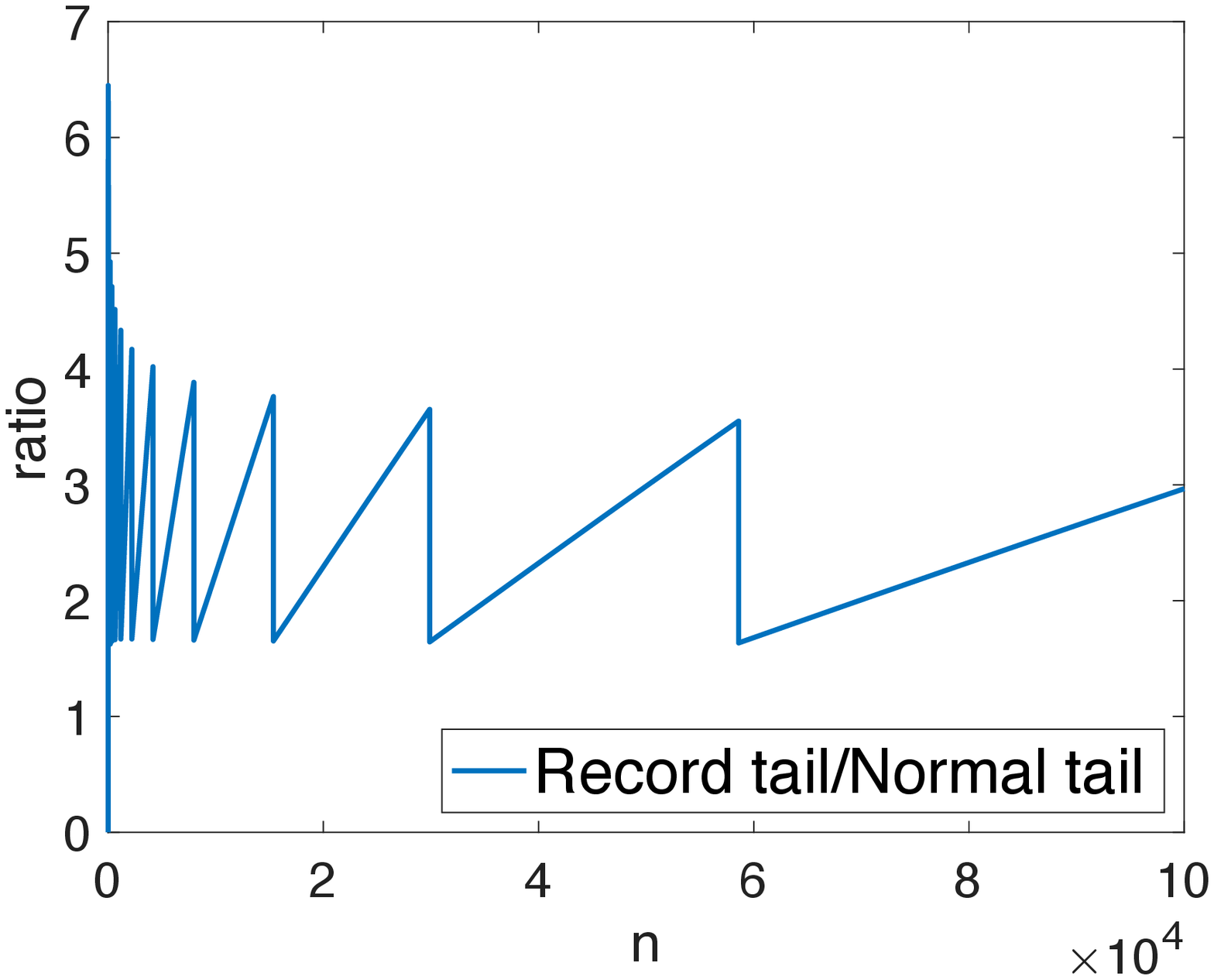}
 \caption{$N(\lambda_n,\sigma_n^2)$ without correction} %
\label{figure2} \end{minipage}
    \end{figure}
    
\begin{figure}[ht]
\begin{minipage}[t]{0.49\linewidth}  % <---
%the four numbers after "trim=" are for left bottom right top
   \includegraphics[trim = 10mm 75mm 5mm 80mm,height=0.26\textheight]{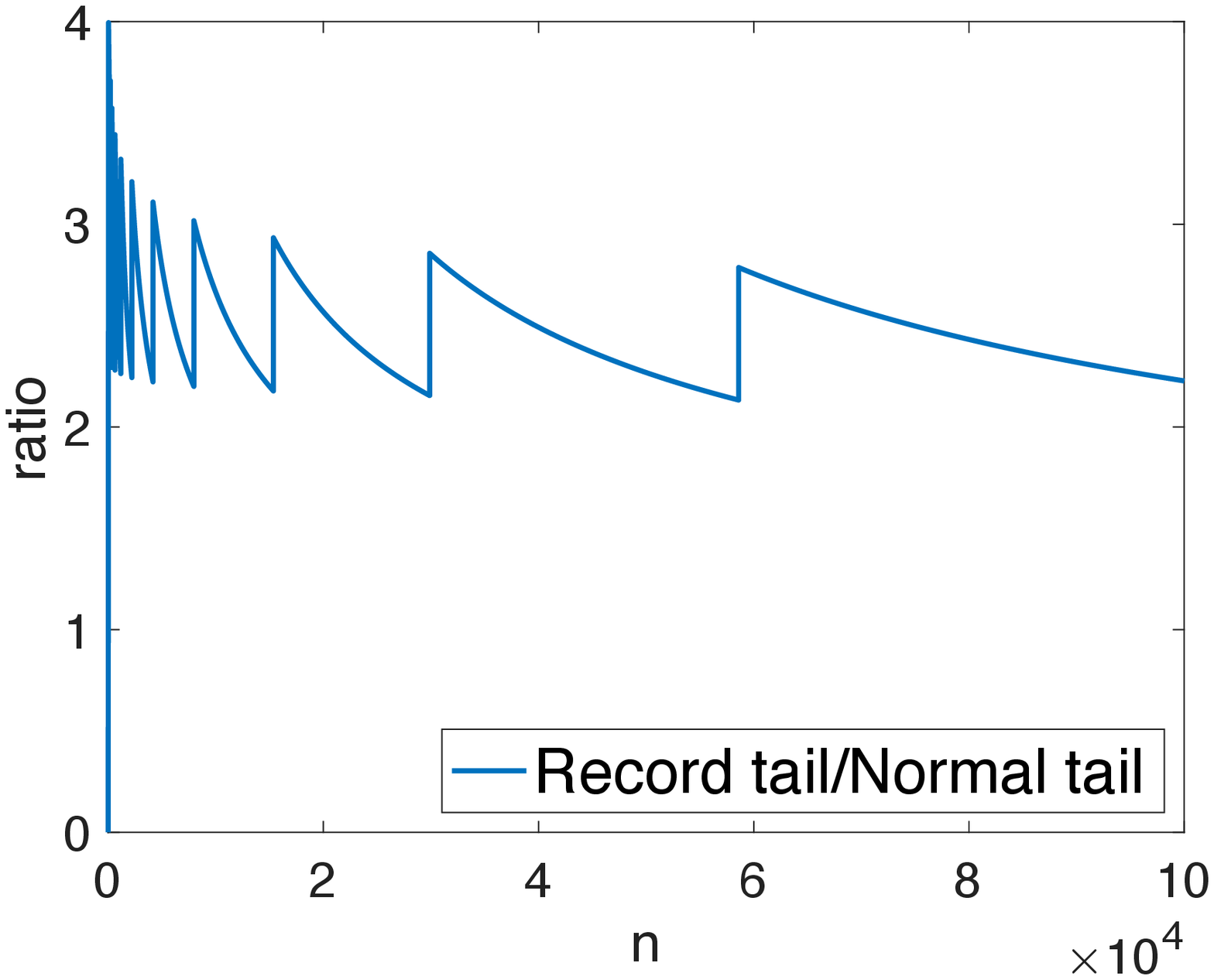} 
\caption{$N(\lambda_n,\sigma_n^2)$ with correction} %
\label{figure3}  \end{minipage}
\hfill
\begin{minipage}[t]{0.5\linewidth}  % <---
%the four numbers after "trim=" are for left bottom right top
   \includegraphics[trim = 10mm 75mm 5mm 80mm,height=0.26\textheight]{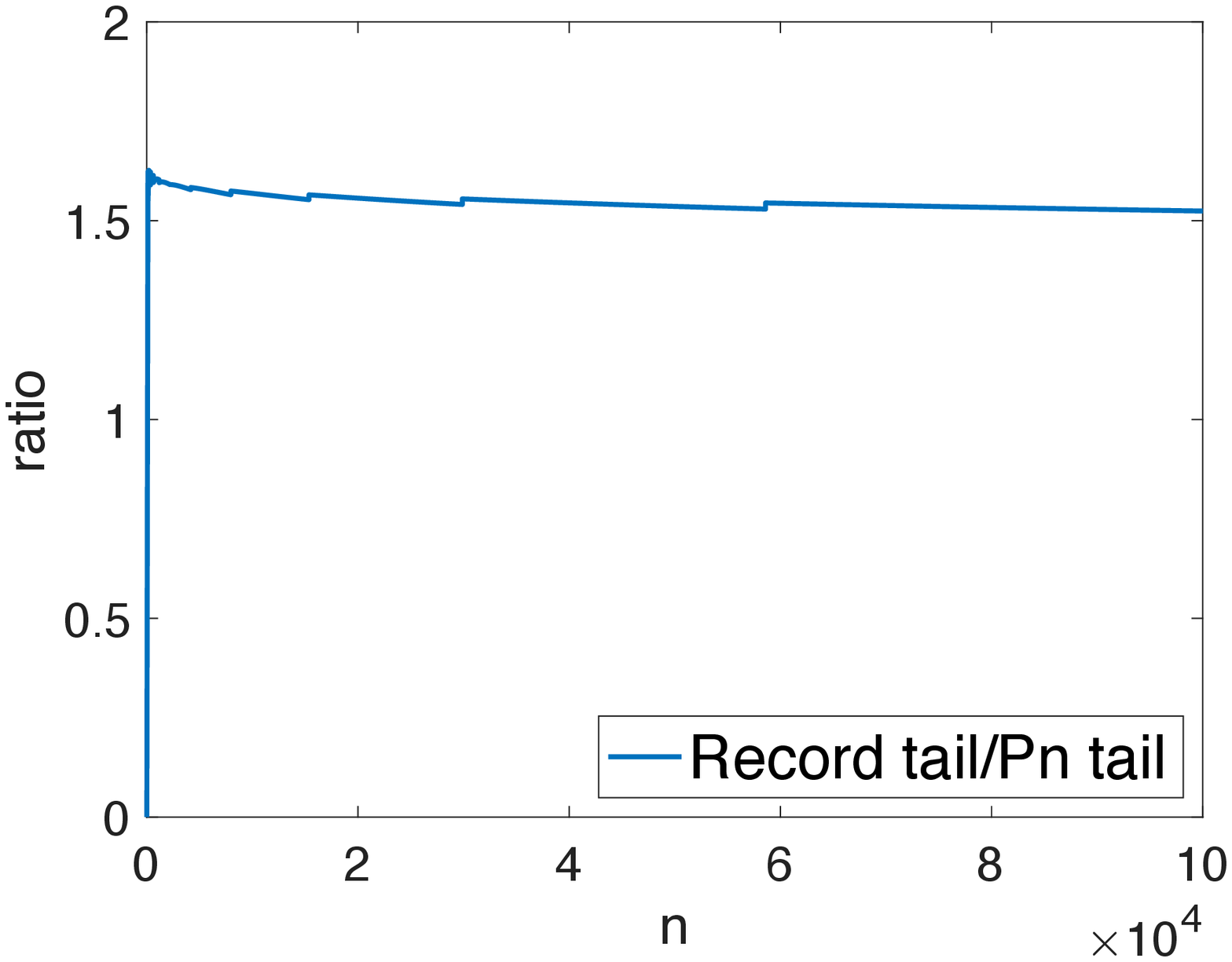}
 \caption{$\Pn(\sigma_n^2)$} %
\label{figure4} \end{minipage}
    \end{figure}
\end{exam}

Example~\ref{exam1} shows that the distribution of the counts of rare events often has a heavier right tail than that of the corresponding normal distribution, 
approximations by the moderate deviations in the normal distribution are generally \red{inferior to those by the moderate deviations in the
Poisson distribution}. The next example says that the parameter of the approximating Poisson distribution suggested in \cite{CC92,BCC95,CFS} is not optimal and some adjustment can significantly improve the quality of approximations by the moderate deviations in the Poisson distribution.

\begin{exam}\label{exam2} {\rm  With $0<p<1$, let $W_n\sim\text{Bi}(n,p)$, $Y_n\sim\Pn(np)$ and $Z\sim N(0,1)$, then for a fixed $x>0$, 
\red{$$\lim_{n\to\infty}\frac{\mathbb{P}(W_n\ge np+x\sqrt{np(1-p)})}{\mathbb{P}(Y_n\ge np+x\sqrt{np(1-p)})}= \frac{\mathbb{P}\left(Z\ge x\right)}{\mathbb{P}(Z\ge x\sqrt{1-p})},$$}
which systematically deviates from 1 as $x$ moves away from 0. The systematic bias can be removed by introducing an adjustment into the approximate models: for a fixed $x>0$, 
\red{$$\lim_{n\to\infty}\frac{\mathbb{P}(W_n\ge np+x\sqrt{np(1-p)})}{\mathbb{P}(Y_n\ge np+x\sqrt{np})}{=}1$$}
or equivalently, with $Y_n'\sim \Pn(np(1-p))$, 
\red{$$\lim_{n\to\infty}\frac{\mathbb{P}(W_n\ge np+x\sqrt{np(1-p)})}{\mathbb{P}(Y_n'\ge np(1-p)+x\sqrt{np(1-p)})}{=}1.$$}
}
\end{exam}

Example~\ref{exam2} suggests that it is more suitable to approximate the right tail probabilities by looking at the number of standard variations away from the mean, which is essentially the original idea of the translated (shifted) Poisson approximation \cite{BX99,Roellin05,Roellin07}. In this paper, we show that {it is indeed better to approximate} the right tail probabilities via the moderate deviations in the {translated} Poisson distribution. 

Our approach does not rely on the boundedness of the Radon-Nikodym derivative as in \cite{CC92,BCC95} or the tacit assumption of well-behaved tail probabilities as in \cite{CFS}, see Remark~\ref{diffofliterature} for more details. For the case of Poisson-binomial, we show in Proposition~\ref{secondpropPnBi} that our approach works for the case that the maximum of the success probabilities of the Bernoulli random variables is not small, such as the distribution of the number of records.

The paper is organised as follows. We state the main results in the context of local dependence, size-biased distribution 
and discrete zero-biased distribution in Section~\ref{secMainresult}. The accuracy of our bounds is illustrated in six examples in Section~\ref{secExamples}. The proofs of the main results are postponed to Section~\ref{secProof} where we also establish Stein's factors for Poisson moderate deviations in Lemma~\ref{lma2}.

\section{The main results}\label{secMainresult}

In this section, we state three theorems on moderate deviations in Poisson approximation, the first is under a local dependent structure, the second is with respect to the size-biased distribution and the last is in terms of the discrete zero-biased distribution. 

We first consider a class of non-negative integer valued random variables $\{X_i:\ i\in\ci\}$ satisfying the local dependent structure (LD2) in \cite{ChenShao04} \red{(see  also \cite{AGG89} for its origin)}. For ease of reading, we quote the definition of (LD2) below. 
\vskip5pt
\begin{tabular}{l} (LD2) For each $i\in\ci$, there exists an $A_i\subset B_i\subset \ci$ such that $X_i$ is independent of \\
\hskip1.2cm$\{X_j:\ j\in A_i^c\}$ and $\{X_i:\ i\in A_i\}$ is independent of $\{X_j:\ j\in B_i^c\}$.
\end{tabular}
\vskip5pt
We set $W=\sum_{i\in\ci} X_i$, $Z_i=\sum_{j\in A_i} X_j$, $Z_i'=\sum_{j\in B_i}X_j$, $W_i=W-Z_i$ and $W_i'=W-Z_i'$. We write 
 \be{\mu_i=\mathbb{E} (X_i),~~\mu=\mathbb{E} (W),~~\sigma^2=\Var{(W)}.} As suggested in Example~\ref{exam2}, we consider $Y\sim\Pn(\lambda)$ approximation to $W-a$ with \red{$\left|\lambda-\sigma^2\right|$ being not too large} and $a=\mu-\lambda$ being an integer so that $k$ in $\mathbb{P}(W-a\ge k)$ and $\mathbb{P}(Y\ge k)$ is in terms of the number of standard deviations of $W$. In principle, the constant $a$ is chosen to minimise the error of approximation, however, our theory is formulated in such a flexible way that other choices of $\lambda$ and $a$ are also acceptable. The three most useful choices of $a$ are $a=0$, $a=\floor {\mu-\sigma^2}$ and $a=\ceil{\mu-\sigma^2}$, where $\floor{\cdot}$ stands for the largest integer in $(-\infty,\cdot]$.

\begin{thm}\label{thm1} With the setup in the preceding paragraph, assume that $\{X_i:\ i\in\ci\}$ satisfies (LD2) and, for each $i$, there exists a $\sigma$-algebra $\cf_i$ such that $\{X_j:\ j\in B_i\}$ is $\cf_i$
measurable. Define
$$\theta_i:=\esssup\max_j\mathbb{P}(W=j|\cf_i),$$
where $\esssup V$ is the essential supremum of the random variable $V$. Then for integer $a<\mu$, $\lambda=\mu-a$ and positive integer $k>\lambda$, we have
\bean{\left\vert\frac{\mathbb{P}(W-a\ge k)}{\mathbb{P}(Y\ge k)}-1\right\vert\le&\bc_2(\lambda,k)\sum_{i\in\ci}\theta_i\left\{{|\mathbb{E}(X_i-\mu_i)Z_i|}\mathbb{E} (Z_i')\right.\non\\
&+\left.\mathbb{E}\left[|X_i-\mu_i|Z_i(Z_i'-Z_i/2-1/2)\right]\right\}\non\\
&+\bc_1(\lambda,k)|\lambda-\sigma^2|+\mathbb{P}(W-a<-1),
\label{esti1}
}
{where, with $F(j)=\mathbb{P}(Y\le j)$, $\overline{F}(j)=\mathbb{P}(Y\ge j)$,}
{\begin{eqnarray}\bc_1(\lambda,k)&:=&\frac{F(k-1)}{k\mathbb{P}(Y=k)}\left\{1-\min\left(\frac{F(k-2)}{F(k-1)}\cdot\frac{\lambda}{k-1},\frac{\overline{F}(k+1)}{\overline{F}(k)}\cdot\frac{k}{\lambda}\right)\right\},\label{Steinconstantnew1}\\
\bc_2(\lambda,k)&:=&\frac{F(k-1)}{k\mathbb{P}(Y=k)}\left(2-\frac{F(k-2)}{F(k-1)}\cdot\frac{\lambda}{k-1}-\frac{\overline{F}(k+1)}{\overline{F}(k)}\cdot\frac{k}{\lambda}\right).\label{Steinconstantnew2}
\end{eqnarray}}
\end{thm}

\begin{rem}\label{remarknaive}
{\rm Both $\bc_1$ and $\bc_2$ can be numerically computed in applications and they can't be generally improved (see the proofs below). 
They are better than the ``naive'' counterparts 
$(1-e^{-\lambda})/(\lambda \mathbb{P}(Y\ge k))$ derived through the total variation bounds in \cite{BE84,BHJ}. Figure~\ref{figure5} provides details of 
$$\mbox{ratio }i:=\bc_i(\lambda,k)/[(1-e^{-\lambda})/(\lambda \mathbb{P}(Y\ge k))], \ i=1,2,$$ for $\lambda=10$,
$k$ from $10$ to $43$. We would like to mention that for large $k$ and/or large $\lambda$, the tail probabilities are so small that the calculation using MATLAB produces unstable results since accumulated computation errors often exceed the tail probabilities\red{, hence more powerful computational tools are needed to achieve the required accuracy or one has to resort to known approximations to the Poisson right tails and point probabilities.}  
\vskip5pt
\begin{figure}[h]
\centering
%\begin{minipage}[t]{0.49\linewidth}  % <---
%the four numbers after "trim=" are for left bottom right top
   \includegraphics[trim = 10mm 55mm 5mm 50mm,height=0.4\textheight]{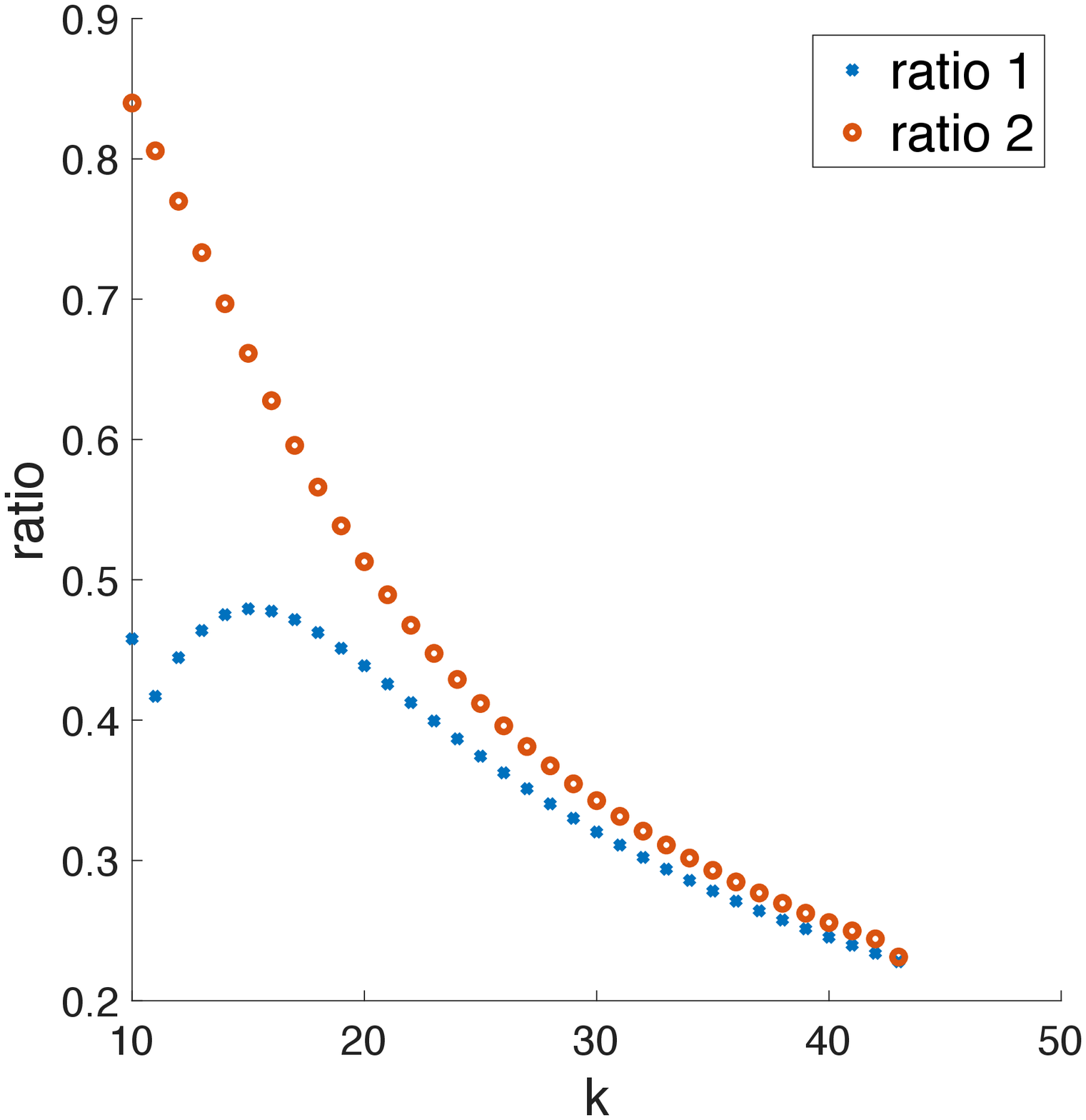}  
\caption{Performance of the bound} %
\label{figure5}  %\end{minipage}
\hfill
%\begin{minipage}[t]{0.5\linewidth}  % <---
%%the four numbers after "trim=" are for left bottom right top
%   \includegraphics[trim = 10mm 75mm 5mm 80mm,height=0.26\textheight]{matching.ps}
% \caption{$C_1(1,k)$} %
%\label{figure6} \end{minipage}
    \end{figure}
}
\end{rem}

\begin{rem}{\rm Due to the discrete nature of Poisson distribution, it seems impossible to analytically simplify $\bc_1$ and $\bc_2$ at negligible costs for the diverse range of $k>\lambda$.
}
\end{rem}

\begin{rem}\label{diffofliterature}{\rm If $\lambda$ is chosen reasonably close to $\sigma^2$ so that $\lambda-\sigma^2$ is bounded, then $\theta_i$ in the bound \Ref{esti1} converges to $0$ when $\sigma^2$ converges to $\infty$. Our bound does not rely on the Radon-Nikodym derivative of $\law(W)$ with respect to $\Pn(\lambda)$, which is the crucial ingredient in   \cite{CC92,BCC95}. On the other hand, the tacit assumption of \cite{CFS} is {that} $\sup_{\lambda\le r\le k}\frac{\mathbb{P}(W\ge r)}{\mathbb{P}(Y\ge r)}$ for $W$ and $Y$ in Theorem~\ref{thm1} is well-behaved and this assumption is hard to verify. The bound \Ref{esti1}, although relatively crude, does not rely on this assumption and covers more general cases.}
\end{rem}

%\begin{rem}
%Assume that \be{\max_i\vert B_i\vert\le m,~~\max_i\vert N_i\vert\le m,} then the upper bound is of order $O(n^{-1/2})$.
%\end{rem}
\begin{cor}\label{cor1}For the sum of independent non-negative integer valued random variables $W=\sum_{i\in\ci}X_i$, let $\theta_i=\max_{j}\mathbb{P}(W-X_i=j)$, $\mu_i=\mathbb{E} X_i$, $\mu=\sum_{i\in\ci} \mu_i$, $\sigma^2=\Var(W)$. For any integer $a<\mu$, let $\lambda=\mu-a$, $Y\sim\Pn(\lambda)$, then for $k> \lambda$,
\bea{&\left\vert\frac{\mathbb{P}(W-a\ge k)}{\mathbb{P}(Y\ge k)}-1\right\vert\non\\
&\le \bc_2(\lambda,k)\sum_{i\in\ci}\theta_i\left\{\mu_i{|\mathbb{E}[X_i(X_i-\mu_i)]|}+\frac{1}{2}\mathbb{E}\left[\vert X_i-\mu_i\vert X_i(X_i-1)\right]\right\}\non\\
&\ \ \ +\bc_1(\lambda,k)|\lambda-\sigma^2|+\mathbb{P}(W-a<-1).%\label{cor1.1}
}\end{cor}

\begin{rem}{\rm We leave $\mathbb{P}(W-a<-1)$ in the upper bound \Ref{esti1}  because the current approach can not remove it from the bound. Nevertheless, it is no more than $1$ and {converges to zero exponentially fast} with suitable choice of $a$. For the sum of independent non-negative integer valued random variables in Corollary~\ref{cor1}, if $a$ is at least less than $\mu$ by a few $\sigma$s, we can use \cite[Theorem~2.7]{CL06} to obtain
\begin{equation}\mathbb{P}(W-a<-1)\le e^{-\frac{(\mu-a+2)^2}{2\sum_{i\in\ci}\mean(X_i^2)}}. \label{cor1-1}
\end{equation}
}
\end{rem}

For any non-negative random variable $W$ with mean $\mu
\in (0,\infty)$ and distribution $dF(w)$, the $W$-{\it size
biased} distribution \cite{Cochran,AG10} is given by
\be{
%\label{size-bias-distribution}
dF^s(w) = \frac{wdF(w)}{\mu}, \quad
\mbox{$w \ge 0$,}}
or equivalently by the characterising equation
$$
\mean [Wg(W)]=\mu \mean g(W^s) \quad \mbox{for all $g$ with
$\mean|Wg(W)|<\infty$.}
$$

\begin{thm}\label{thm1.2}Let $W$ be a non-negative integer-valued random variable with mean $\mu$ and variance $\sigma^2$, $a<\mu$ be an integer, $\lambda=\mu-a$. Then for integer $k>\lambda$, we have \bean{
\left\vert\frac{\mathbb{P}(W-a\ge k)}{\mathbb{P}(Y\ge k)}-1\right\vert\le&\bc_1(\lambda,k)\left\{\mu\mathbb{E}\vert W+1-W^s\vert+\vert\mu-\lambda\vert\right\}\non\\
&+\mathbb{P}(W-a<-1),\label{estimate1.2}
}where $Y\sim\Pn(\lambda)$.
\end{thm}

\begin{rem}{\rm  Theorem~\ref{thm1.2} improves \cite[Theorem~3]{CFS} in a number of ways, with less restrictive conditions %, sharper estimates 
and no unspecified constants.
}
\end{rem}

The next theorem is based on the discrete zero-biased distribution defined in \cite{Goldstein06} and the approach is very similar to that in \cite{CFSb}. For an integer valued random variable $V$ with mean $\mu$ and finite variance $\sigma^2$, we say that
$V^\star$ has the discrete $V$-zero biased distribution \cite[Definition~2.1]{Goldstein06} if, for
all bounded functions $g:\ \Z:=\{0,\pm1,\pm2,\dots\} \rightarrow \R$ with $\mathbb{E}|Vg(V)|<\infty$,
$$\mean [(V-\mu)g(V)] =
\sigma^2\mean\Delta g(V^\star),%\label{star}
$$
{where $\Delta f(i):=f(i+1)-f(i)$.} 
 
\begin{thm}\label{thm2}Let $W$ be a non-negative integer-valued random variable with mean $\mu$, variance $\sigma^2$, $a<\mu$ be an integer,  and $W^\star$ have the discrete $W$-zero biased distribution and be defined on the same probability space as $W$.
Set $R=W^\star-W$ and define
$$\theta_R=\max_{j}\mathbb{P}(W=j|R).$$
Then, for integer $k>\lambda$, with $\lambda=\mu-a>0$, we have 
\bean{
&\left\vert\frac{\mathbb{P}(W-a\ge k)}{\mathbb{P}(Y\ge k)}-1\right\vert\non\\
&\le\bc_2(\lambda,k)\sigma^2\mathbb{E}[|R|\theta_R]+\bc_1(\lambda,k)|\lambda-\sigma^2|\lambda^{-1}+\mathbb{P}(W-a<-1),\label{estimate2}
}where $Y\sim\Pn(\lambda)$.
\end{thm}

\section{Examples}\label{secExamples}

As many applications of Poisson approximation rely on size biased distributions, we begin with a review of some facts about size biasing. 

Size biasing has been of considerable interest for many decades (see \cite{BHJ}, \cite{Ross11}, \cite{AGK13} and references therein). In the context of the sum of Bernoulli random variables, 
its size biasing is particularly simple. More precisely, 
if $\{X_i:\ i\in\ci\}$ is a family of Bernoulli random variables with $\mathbb{P}(X_i=1)=p_i$, then the size biased distribution of $W=\sum_{i\in\ci}X_i$
is 
\ben{W^s=\sum_{j\ne I}X_j^{(I)}+1,\label{sizebiasing}}
where 
$$\law(\{X_j^{(i)}:\ j\in\ci\})=\law(\{X_j:\ j\in\ci\}\vert X_i=1),$$
$I$ is a random element independent of $\{\{X_j^{(i)}:\ j\in\ci\}:\ i\in\ci\}$ having distribution $\mathbb{P}(I=i)=\frac{p_i}{\mathbb{E} W}$, $i\in\ci$.
Moreover, $\{X_i:\ i\in\ci\}$ are said to be negatively related (resp. positively related) \cite[p.~24]{BHJ} if one can construct $\{\{X_j^{(i)}:\ j\in\ci\}:\ i\in\ci\}$ such that $X_j^{(i)}\le$ (resp. $\ge$) $X_j$ for all $j\ne i$.  When $\{X_i:\ i\in\ci\}$ are negatively related, we have
\bean{\mathbb{E}\vert W+1-W^s\vert
=\mathbb{E}(W+1-W^s)=\mu^{-1}(\mu-\sigma^2)\label{monotonecoupling1},}
where $\mu=\mathbb{E} W$ and $\sigma^2=\var(W)$.
On the other hand, if $\{X_i:\ i\in\ci\}$ are positively related, then 
\bean{\mathbb{E}\vert W+1-W^s\vert&=\mathbb{E} \left\vert\sum_{j\ne I}(X_j^{(I)}-X_j)-X_I\right\vert\non\\
&\le\mathbb{E}\left\{\sum_{j\ne I}(X_j^{(I)}-X_j)+X_I\right\}\non\\
&=\mathbb{E}(W^s-W-1)+2\mu^{-1}\sum_{i\in\ci}p^2_i\non\\
&=\mu^{-1}(\sigma^2-\mu)+2\mu^{-1}\sum_{i\in\ci}p^2_i.\label{monotonecoupling2}}

\subsection{Poisson-binomial trials}

Let $\{X_i,~1\le i\le n\}$ be independent Bernoulli random variables with $\mathbb{P}(X_i=1)=p_i\in (0,1)$, $W=\sum_{i=1}^nX_i$, {$\mu=\mathbb{E} W$} and {$\mu_2=\sum_{i=1}^np_i^2$}. When $\tilde p:=\max_{1\le i\le n}p_i\to 0$, the large deviation of $W$ is investigated in \cite{CC92,BCC95} with precise asymptotic order. We give two results for this particular case without the assumption $\tilde p$ being small, the first is direct consequences of the general results in Section~\ref{secMainresult} and the second is based on our approach using a more fine-tuned analysis and well-studied properties of the tail behaviour of $W$.

\begin{prop} \red{Recalling $\bc_1$ and $\bc_2$ in \Ref{Steinconstantnew1} and \Ref{Steinconstantnew2}, for any integer $k>\mu$, we have}
\ben{\left\vert\frac{\mathbb{P}(W\ge k)}{\Pn(\mu)([k,\infty))}-1\right\vert\le\bc_1(\mu,k)\mu_2\label{Poissonbinomial0}}%\Comment{Please check the details}
and, with {$a=\floor{\mu_2}$ and $\lambda:=\mu-a$},
\bean{\left\vert\frac{\mathbb{P}(W-a\ge k)}{\Pn(\lambda)([k,\infty))}-1\right\vert\le&
\frac{\bc_2(\lambda,k)\sum_{i=1}^np_i^2(1-p_i)}{1\vee \sqrt{(\sum_{i=1}^np_i\wedge(1-p_i)-1/4)\pi/2}}\non\\
&+\bc_1(\lambda,k)|\lambda-\sigma^2|+{e^{-(\lambda+2)^2/(2\mu)}}.\label{PoissonBino}} 
\end{prop}

\noindent{\bf Proof} The claim \Ref{Poissonbinomial0} is a consequence of Theorem~\ref{thm1.2} with $a=0$ and $\mu\mathbb{E} |W+1-W^s|=\sum_{i=1}^np_i^2$, as shown in \Ref{monotonecoupling1}.

 The bound \Ref{PoissonBino} is a special case of Corollary~\ref{cor1}. Since $\law(W_i)$ is unimodal, \cite[Corollary~1.6]{MR07} says that
\bean{\theta_i= \dtv(W_i,W_i+1)&\le1\wedge\left\{\sqrt{\frac{2}{\pi}}\left(\frac{1}{4}+\sum_{j\ne i}p_j\wedge(1-p_i)\right)^{-1/2}\right\}\non\\
&\le 1\wedge \left\{\sqrt{\frac{2}{\pi}}\left(\sum_{i=1}^np_i\wedge(1-p_i)-1/4\right)^{-1/2}\right\}.
\label{MRmaximumprob}}
On the other hand, $\mathbb{E}(X_i^2)=p_i$, hence the upper bound \Ref{PoissonBino} is an immediate consequence of Corollary~\ref{cor1} and \Ref{cor1-1}. 

One can also use Theorem~\ref{thm2} to obtain the same bound. More precisely, according to the construction of the discrete zero-biased distribution suggested in \cite{Goldstein06}, let $I$ be a random variable independent of 
$\{X_i,~1\le i\le n\}$ with distribution $\mathbb{P}(I=i)=p_i(1-p_i)/\sigma^2$ for $1\le i\le n$, then we can write $W^\star=W-X_I$, giving $R=-X_I$. We then apply \Ref{MRmaximumprob} to bound $\theta_R$ as 
$$\theta_R=\max_{j}\mathbb{P}(W=j|R)\le \sqrt{\frac{2}{\pi}}\left(\sum_{i=1}^np_i\wedge(1-p_i)-1/4\right)^{-1/2},$$
and a routine calculation gives $\mathbb{E}|R|=\sum_{i=1}^np_i^2(1-p_i)/\sigma^2$, hence \Ref{PoissonBino} follows from \Ref{estimate2} and \Ref{cor1-1}. \qed
\vskip10pt
\begin{prop}\label{secondpropPnBi} Define$$M:=M(p_1,\dots,p_n)=\left\{\begin{array}{ll}
e^\mu,&\mbox{ if }0<\mu<1,\\
e^{13/12}\sqrt{2\pi}\left(1-\mu_2/\mu\right)^{-1/2},&\mbox{ if }\mu\ge 1,
\end{array}\right.$$
then for any integer $k$ with $x:=(k-\mu)/\sqrt{\mu}\ge 1$, we have
\ben{0>\frac{\mathbb{P}(W\ge k)}{\Pn(\mu)([k,\infty))}-1>- 2M(\mu_2/\mu)\left(x^2+{1}+4x{\sqrt{\frac{1-e^{-\mu}}\mu}}\right).\label{Poissonbinomial1}}
\end{prop}

The proof relies on more information of the solutions of Stein's equation and it is postponed to the end of Section~\ref{secProof}. The bound \Ref{Poissonbinomial1} improves \cite[(3.1)]{CFS} \red{in two aspects: it contains no unspecified constants and it does not require $\tilde p$ being small}. For the distribution of the number $S_n$ of records, the large deviation results in \cite{BCC95} do not apply. However, recalling that $\lambda_n=\sum_{i=2}^n\frac1i$, we apply Proposition~\ref{secondpropPnBi} with the harmonic series $\lambda_n=\sum_{i=2}^n\frac1i\ge \ln n+\gamma-1$ and the Riemann zeta function $\sum_{i=2}^n\frac1{i^2}\le \sum_{i=2}^\infty\frac1{i^2}=\frac{\pi^2}6-1$ to get the following estimate.

\begin{cor}  For any integer $k$ with $x:=(k-\lambda_n)/\sqrt{\lambda_n}\ge 1$, we have
$$0>\frac{\mathbb{P}(S_n\ge k)}{\Pn(\lambda_n)([k,\infty))}-1> -\frac{2e^{13/12}\sqrt{2\pi}(\pi^2/6-1)}{\sqrt{(\ln n+\gamma-1)(\ln n+\gamma-\pi^2/6)}}\left(x^2+{1}+\frac{4x}{\sqrt{\ln n+\gamma-1}}\right),$$ 
where $\gamma$ is Euler's constant.
\end{cor}

{\begin{rem}{\rm We conjecture that, with $a=\floor{\mu_2}$ and $\lambda:=\mu-a$, the bound in \Ref{PoissonBino}
can be significantly improved and the better estimate is likely dependent on the Radon-Nikodym derivative bound $\sup_{r\ge 0}\frac{\mathbb{P}(W-a=r)}{\Pn(\lambda)(\{r\})}$. }
\end{rem}}
 
\subsection{Matching problem}%\label{matchinsection}

For a fixed $n$, let $\pi$ be a uniform random permutation of $\{1,\dots,n\}$, $W=\sjn\bone_{\{j=\pi(j)\}}$ be the number of fixed points in the permutation. 

\begin{prop}\label{propmatching} For the random variable $W$ defined above and any integer $k\ge 2$, we have
\ben{\left\vert\frac{\mathbb{P}(W\ge k)}{\Pn(1)([k,\infty))}-1\right\vert\le{\frac{2}{n}\bc_1(1,k)}.\label{matching1-1}}
%{See figure~\ref{figure6} for the performance of $\bc_1(1,k)$ for $1\le k\le 6$.}
\end{prop}

\noindent{\bf Proof of Proposition~\ref{propmatching}} In this case, the size-biased distribution $\law(W^s)$ can be coupled with $W$ as follows \cite{CDM}. Let $I$ be uniformly distributed on $\{1,2,\dots,n\}$, independent of $\pi$, and define \be{\pi^s(j)=\begin{cases}
I&\text{if~}j=I,\\
\pi(I)&\text{if~}j=\pi^{-1}(I),\\
\pi(j)&\text{otherwise}.\end{cases}}
Set $W^s=\sjn\bone_{\{j=\pi^s(j)\}}$, one can easily verify that $W^s$ has the size-biased distribution of $W$. Also, we can check that $\mathbb{E} W=\Var(W)=1$, giving $\mathbb{E} W^s=2$. 
Let $\Delta=W+1-W^s$, using the above construction of $W^s$, we can conclude that $\Delta$ takes values in $\{-1,0,1\}$ and  
$\mathbb{P}(\Delta=1\vert W)=W/n$. Since $\mathbb{E}\Delta=0,$ we have $\mathbb{P}(\Delta=1)=\mathbb{P}(\Delta=-1)$, and $\mathbb{E}\vert\Delta\vert=2/n$. On the other hand, $\lambda=\mu$ allows us to get rid of the second term in \Ref{estimate1.2}. By Theorem \ref{thm1.2} with $a=0$, $\lambda=\mu=1$, the claim follows. \qed

\begin{rem}{\rm The bound \Ref{matching1-1} contains no unknown constants and improves the bound of \cite[\S3.3]{CFS}. \ignore{If we derive a naive bound using  
$$\dtv(\law(W),\Pn(1))\le 2(1-e^{-1})/n$$ 
(see \cite{CDM}), then $\bc_1(1,k)$ is replaced with $\bc_1'(1,k)=(1-e^{-1})/{\Pn(1)([k,\infty))}$. We list the differences between these two constants for $k=1,\dots,9$ in Table~\ref{table1}.
\vskip10pt
\begingroup
\setlength{\tabcolsep}{4.6pt} % Default value: 6pt
\renewcommand{\arraystretch}{1.5} % Default value: 1
\centerline{\begin{tabular}{|c|ccccccccc|}
\hline
$k$            &1 &2 &3 &4 &5 &6 &7 &8 &9\\
\hline
$\bc_1(1,k)$   &0.58   &1.00      &3.00 &11.00   &49.00     &261.00   &1631.00       &11743.00        &95901.00\\
$\bc_1'(1,k)$  &1.00   &2.39      &7.87 &33.29   &172.72   &1063.85 &7593.85   &61675.13   &561783.77\\
\hline 
$\frac{\bc_1'(1,k)}{\bc_1(1,k)}$&1.72&2.39&2.62&3.03& 3.52&4.08&4.66&5.25&5.86\\ [2pt]
\hline
\end{tabular}}
\captionof{table}{Constants for matching problem}\label{table1}
\endgroup
\vskip10pt}
}
\end{rem}

\subsection{Occupancy Problem}
The occupancy problem has a long history dating back to the early development of probability theory. General references on this subject can be found in classics, e.g., \cite[Vol~1, Chapter~2]{Feller68} and \cite[Chapter 6]{BHJ}. 

The occupancy problem can be formulated as follows. Let $l$ balls be thrown independently of each other into $n$ boxes uniformly. Let $X_i$ be the indicator variable of the event that $i$-th box being empty, so the number of empty boxes can be written as $W=\sum_{i=1}^nX_i$. Noting that $p:=\mathbb{E} X_i=\left(1-\frac{1}{n}\right)^l$, direct computation gives \be{\mu:=\mathbb{E} W=np,}
\bea{\sigma^2:=\Var(W)=\mu-\mu^2+\mu(n-1)\left(1-\frac{1}{n-1}\right)^l.}

\begin{prop}\label{propoccupancy}
For the random variable $W$ defined above and any integer $k>\mu$, we have
\ben{\left\vert\frac{\mathbb{P}(W\ge k)}{\mathbb{P}(Y\ge k)}-1\right\vert\le\bc_1(\mu,k)\mu\left[\mu-(n-1)\left(1-\frac{1}{n-1}\right)^l\right],
\label{occupancyadd1}} where $Y\sim \Pn(\mu)$.
\end{prop}
\noindent{\bf Proof of Proposition~\ref{propoccupancy}} For the sake of completeness, we provide the following proof which is essentially a repeat of \cite[p.~23]{BHJ}.
From the construction of $W$-size biased distribution in \Ref{sizebiasing}, we can construct a coupling as follows. Let $I$ be uniform on $\{1,\dots,n\}$, that is, we randomly pick one box with equal probability. If the selected box is not 
empty, we redistribute all balls in the box randomly into the other $n-1$ boxes with equal probability $1/(n-1)$. Define $X_{j}^{(i)}$ as the indicator of the event that the box being selected is $i$, and after the redistribution, box $j$ is empty. With this coupling in mind, one can verify that $\{X_i\}$ is negatively related so it follows from 
\Ref{monotonecoupling1} that
\be{
\mathbb{E}\vert W+1-W^s\vert=\mu-(n-1)\left(1-\frac{1}{n-1}\right)^l.
} 
{Now, applying Theorem~\ref{thm1.2} with $a=0$ yields} \Ref{occupancyadd1}. \qed

\subsection{Birthday problem}

The classical birthday problem is essentially a variant of the occupancy problem. For this reason, we throw $l$ balls independently and equally likely into $n$ boxes and let $X_{ij}$ be the indicator random variable of the event that ball $i$ and ball $j$ fall into the same box. The number of pairs of balls going into the same boxes (i.e., the number of pairs of people having the same birthdays) can be written as $W=\sum_{i<j}X_{ij}$. Define $p=\mathbb{E} X_{ij}=\frac{1}{n}$, so $\mu=\mathbb{E} W={l\choose 2}p.$ \cite{CDM} \red{give the following construction of $W^s$: label the balls from $1$ to $l$, randomly choose two balls $J_1$ and $J_2$ and move ball $J_1$ into the box that $J_2$ is in, then $W$ is the number of pairs of balls before the move while $W^s$ is the number of pairs of balls after the move. Let $E$ be the event that $J_1$ and $J_2$ are from the same box. When $E$ occurs, $W^s=W$ so $\vert W+1-W^s\vert=1$; otherwise, $J_1$ and $J_2$ are from different boxes with $B_1$ and $B_2$ balls respectively, giving
\be{ W+1-W^s=B_1-B_2.}
Hence,
\bea{\mathbb{E}\vert W+1-W^s\vert&=\mathbb{P}(E)+\mathbb{E}[\vert W+1-W^s\vert\vert E^c]\mathbb{P}(E^c)\\
&\le\frac1n+\mathbb{E}\vert B_1-B_2\vert\\
&\le\frac1n+\mathbb{E} (B_1+B_2)=\frac{1+2l}{n}.}  
This, together with Theorem~\ref{thm1.2} and $a=0$, gives the following Proposition.}

\begin{prop}
For the random variable $W$ defined above and any integer $k>\mu$, we have
\be{\left\vert\frac{\mathbb{P}(W\ge k)}{\mathbb{P}(Y\ge k)}-1\right\vert\le\bc_1(\mu,k)\mu\frac{1+2l}{n},} where $Y\sim \Pn(\mu)$.
\end{prop}

\ignore{\begin{rem}
If $l=c\sqrt{n}$ for $c$ fixed and $n$ large, from Theorem~\ref{thm1.2} we have an estimate of the order $O(n^{-1/2})$. 
\end{rem}}

\subsection{\red{Triangles in the \ER~random graph}}
Let $G=G(n,p)$ be an \ER~random graph on $n$ vertices with edge probability $p$. Let $K_n$ be the complete graph on $n$ vertices, and $\Gamma$ be the set of all triangles in $K_n$. For $\a\in\Gamma$, let $X_\a$ be the indicator that there is a triangle in $G$ at $\a$, i.e. \be{X_\a=\bone_{\{\a\subset G\}}.} Therefore the number of triangles in $G$ can be represented as $W=\sum_{\a\in\Gamma}X_\a$. \red{Clearly, $X_\a$ is independent of $X_\b$ if $\a$ and $\b$ don't share a common edge. By analysing} the numbers of shared edges, we obtain (see \cite[p.~255]{Ross11})\be{\mu=\mathbb{E} W={n\choose 3}p^3,}\be{\sigma^2=\Var(W)={n \choose 3}p^3[1-p^3+3(n-3)(p^2-p^3)].}

\begin{prop}\label{proprandomgraph}
For the random variable $W$ defined above and any integer $k>\mu$, we have
\ben{\left\vert\frac{\mathbb{P}(W\ge k)}{\mathbb{P}(Y\ge k)}-1\right\vert\le\bc_1(\mu,k)\mu\left(3(n-3)p^2(1-p)+p^3\right),\label{randomgraph}} where $Y\sim \Pn(\mu)$.
\end{prop}
\noindent{\bf Proof of Proposition~\ref{proprandomgraph}} The following proof is a special version of the general argument in \cite[p.~89]{BHJ}. \red{Since $X_\a$ and $X_\b$ are independent if $\a$ and $\b$ have no common edges,} a size biased distribution of $W$ can be constructed as follows. Let \be{X_\b^{(\a)}:=\bone_{\{\b\subset G\cup\a\}},\ \b\in\Gamma,}
then $\law(\{X_\b^{(\a)},\b\ne\a\})=\law(\{X_\b,\b\ne\a\}|X_\a=1)$. Here the union of graphs is in the sense of set operation of their vertices and edges. Let $I$ be a random element taking values in $\Gamma$ with equal probability and be independent of $\law(\{X_\b^{(\a)},\a,\b\})$, then we can write $W^s=\sum_{\b\ne I}X_\b^{(I)}+1$.  Because $X_\b^{(\a)}\ge X_\b$ {for all $\b\in\Gamma$}, \Ref{monotonecoupling2} implies
\bea{\mathbb{E}\vert W+1-W^s\vert&\le
%E \left\vert\sum_{\b\ne I}(X_\b^{(I)}-X_\b)-X_I\right\vert\\
%&\le\mu^{-1}\sum_{\a\in\Gamma}p^3\mathbb{E}(W_\a-1-W+2X_\a)\\
%\mathbb{E}(W^s-W-1)+2p^3\\
\mu^{-1}(\sigma^2-\mu+2\mu p^3)\\
&=3(n-3)p^2(1-p)+p^3.}
The claim follows from Theorem~\ref{thm1.2} with $a=0$. \qed

\begin{rem}
Since $\mu={n\choose 3}p^3$, if $p=O(1/n)$, then the error bound \Ref{randomgraph} is of the same order $O(1/n)$.
\end{rem}

\subsection{2-runs}%\Comment{needs work in this section}

Let $\{\xi_i,\dots,\xi_n\}$ be i.i.d. $Bernoulli(p)$ random variables with $n\ge 9$, $p<2/3$. For each $1\le i\le n$, define $X_i=\xi_i\xi_{i+1}$ and, to avoid edge effects, we define $\xi_{j+n}=\xi_j$ for $-3\le j\le n$. The number of $2$-runs in the Bernoulli sequence is defined as $W=\sum_{i=1}^nX_i$, then $\mu=np^2$ and variance $\sigma^2=np^2(1-p)(3p+1)$. 
{\begin{prop} 
For any integer $k>\mu$,
\ben{\left\vert\frac{\mathbb{P}(W_n\ge k)}{\Pn(\mu)([k,\infty))}-1\right\vert\le \bc_1(\mu,k)np^3(2-p).\label{2runsadd1}}
With $a:=\floor{np^3(3p-2)}$, $\lambda=\mu-a$, then for any integer $k>\lambda$,
%\bean{&\left\vert\frac{\mathbb{P}(W_n-a\ge k)}{\mathbb{P}(Y\ge k)}-1\right\vert\non\\
%&\le{
%\frac{9.2(4xe^{x^2+1}+1)(1+5p)}{(1+2p-3p^2)\sqrt{(n-8)(1-p)^3}}+\frac{3xe^{x^2+1}}{1\vee[np^2(1+2p-3p^2)]}},\label{2runs2}
%}where $x:=\frac{k-\lambda}{\sqrt{\lambda}}$
\ben{\left\vert\frac{\mathbb{P}(W_n-a\ge k)}{\Pn(\lambda)([k,\infty))}-1\right\vert\le\bc_2(\lambda,k)\frac{9.2np^2(1+5p)}{\sqrt{(n-8)(1-p)^3}}+\bc_1(\lambda,k)(1\wedge\lambda).\label{2runs2}}
\end{prop}}

\noindent{\bf Proof} For \Ref{2runsadd1}, we apply Theorem~\ref{thm1.2} with $a=0$, {$$X_j^{(i)}=\left\{\begin{array}{ll}
X_j,&\mbox{ if }|j-i|\ge 2,\\
\xi_j,&\mbox{ if }j=i-1,\\
\xi_{j+1},&\mbox{ if }j=i+1,\\
1,&\mbox{ if }j=i,
\end{array}\right.$$
$I$ a uniform random variable on $\{1,\dots,n\}$ independent of $\{X_j^{(i)}\}$,}
and  
$$W^s=\sum_{j\ne I}X_j^{(I)}+1,$$ giving 
\bea{\mathbb{E}\vert W+1-W^s\vert&=\mathbb{E}\vert X_{I-1}+X_I+X_{I+1}-\xi_{I-1}-\xi_{I+2}\vert\non\\
&=\mathbb{E}\vert\xi_{i-1}\xi_{i}+\xi_{i}\xi_{i+1}+\xi_{i+1}\xi_{i+2}-\xi_{i-1}-\xi_{i+2}\vert\non\\
&=p(2-p).}

Apropos of \Ref{2runs2}, 
we make use of Theorem~\ref{thm1}. To this end, let
$A_i=\{i-1,i,i+1\}$, $B_i=\{i-2,i-1,i,i+1,i+2\}$, $\cf_i=\sigma\{\xi_j:\ i-2\le j\le i+3\}$, then 
\cite[Lemma~5.1]{BX99} with $\alpha_j=0$ or $1$ for $j=i-2,\cdots,i+5$ gives
\bea{\theta_i&\le\dtv\left(W,W+1|\cf_i\right)\le\frac{2.3}{\sqrt{(n-8)p^2(1-p)^3}}.} 
On the other hand, $\mathbb{E}(Z_i')=5p^2$, {$|\mathbb{E}((X_i-\mu_i)Z_i)|\le \mathbb{E}(Z_i)=3p^2$,} 
$$\mathbb{E}[|X_i-\mu_i|Z_i(Z_i'-Z_i/2-1/2)]\le \mathbb{E}[Z_i(Z_i'-Z_i/2-1/2)]=4p^3+5p^4,$$
and $|\lambda-\sigma^2|\lambda^{-1}\le 1\wedge (\lambda^{-1})$, $a=\floor{np^3(3p-2)}\le0$, {$\lambda\ge \sigma^2$,} hence $\mathbb{P}(W-a<-1)=0$
and \Ref{2runs2} follows from Theorem~\ref{thm1} by collecting these terms. 
\qed

\ignore{\begin{rem}{\rm \cite{BX99} use compound Poisson signed measures with two parameters to approximate the distribution of $2$-runs with an error of order $O(n^{-1/2})$, and the topic is  also studied in \cite{Roellin05} by using a translated Poisson approximation, giving the same order of approximation error. The bound \Ref{2runs2} of moderate deviation approximation is again 
of the same order but with an improvement. At the cost of more complexity, the proof can be extended to study $k$-runs for $k\ge 3$ with unequal probabilities in the Bernoulli trials.}
\end{rem}}

\section{The proofs of the main results}\label{secProof}

The celebrated Stein-Chen method \cite{Chen75} is based on the observation that a non-negative random variable $Y\sim\Pn(\lambda)$
if and only if $\mathbb{E}[\lambda f(Y+1)-Yf(Y)]=0$ for all bounded functions $f:~\Z_+:=\{0,1,2,\dots\}\rightarrow\R$, leading to a Stein identity for Poisson approximation as
\ben{\label{steinid}\lambda f(j+1)-jf(j)=h(j)-\Pn(\lambda)\{h\},~~j\ge 0,}
where $\Pn(\lambda)\{h\}:=\mathbb{E} h(Y)$. Since $f(0)$ plays no role in Stein's equation, we set $f(0)=f(1)$ and $f(j)=0$ for $j<0$. The following Lemma plays the key role in the proofs of the main results and it enables us to circumvent checking the moment condition \Ref{momentgeneratingcondition} which seems to be inevitable in the existing procedure for proving moderate deviation theorems.

\begin{lma}\label{lma2} For fixed {$k\in\Z_+$}, let $h=\bone_{[k,\infty)}$. With {$\pi_\cdot=\Pn(\lambda)(\{\cdot\})$,} $\Delta f(i)=f(i+1)-f(i)$ and $\Delta^2f=\Delta(\Delta f)$, the solution $f:=f_h$ of the Stein equation \Ref{steinid} has the following properties:

{\begin{description}
\item{(i)} $\|f\|:=\sup_{i\in\Z_+}\vert f(i)\vert=\bc_{0}(\lambda,k)\Pn(\lambda)\{h\}$, where $\bc_{0}(\lambda,k):=\frac{F(k-1)}{k\pi_k}$;
\item{(ii)} $\Delta f(i)$ is negative and decreasing in $i\le k-1$; and positive and decreasing in $i\ge k$;
\item{(iii)} $\|\Delta f\|_{k-}:=\sup_{i\le k-1}\vert\Delta f(i)\vert =\bc_{1-}(\lambda,k)\Pn(\lambda)\{h\}$ and
$\|\Delta f\|_{k+}:=\sup_{i\ge k}\vert\Delta f(i)\vert= \bc_{1+}(\lambda,k)\Pn(\lambda)\{h\}$, where 
$$\bc_{1-}(\lambda,k):=\frac{F(k-1)}{k\pi_k}\left(1-\frac{F(k-2)}{F(k-1)}\cdot\frac{\lambda}{k-1}\right)$$
 and  $$\bc_{1+}(\lambda,k):=\frac{F(k-1)}{k\pi_k}\left(1-\frac{\overline{F}(k+1)}{\overline{F}(k)}\cdot\frac{k}{\lambda}\right);$$
\item{(iv)} $\|\Delta f\|:=\sup_{i\in\Z_+}\vert\Delta f(i)\vert=\bc_1(\lambda,k)\Pn(\lambda)\{h\}$ and $\|\Delta^2 f\|:=\sup_{i\in\Z_+}\vert\Delta^2 f(i)\vert=\bc_2(\lambda,k) \Pn(\lambda)\{h\}$;
\end{description}}
\red{where $\bc_1$ and $\bc_2$ are defined in \Ref{Steinconstantnew1} and \Ref{Steinconstantnew2}.}
\end{lma}

{For $k>\l$, death rates are bigger than the birth rate, so it seems intuitively obvious that $\tau_{k}^-$ is stochastically less than or equal to $\tau_{k-2}^+$ for such $k$.  In view of representation \Ref{conj} and $f(k)<0$ as shown in \Ref{sol3}, this is equivalent to $\bc_{1-}(\lambda,k)> \bc_{1+}(\lambda,k)$, leading to the following conjecture.
\begin{con}{\rm We conjecture that $\bc_{1-}(\lambda,k)> \bc_{1+}(\lambda,k)$ for all $k>\lambda$ and the gap increases exponentially as a function of $k-\lambda$.}
\end{con}}

\noindent{\bf Proof of Lemma~\ref{lma2}} We build our argument on the birth-death process representation of the solution
\ben{\label{sol1}
f(i)=-\int_0^{\infty}\mathbb{E}\left[h(Z_i(t))-h(Z_{i-1}(t))\right]dt,\mbox{ for }i\ge 1,
} where $Z_n(t)$ is a birth-death process with birth rate $\lambda$, unit per capita death rate and initial state $Z_n(0)=n$
\cite{Barbour88,BB92,BX01}. For convenience, we adopt the notation in \cite{BX01}: for $i,j\in\Z_+$, define\be{\tau_{ij}=\inf\{t:Z_i(t)=j\},~~\tau^+_j=\tau_{j,j+1},~~\tau^-_j=\tau_{j,j-1}}and \be{\tpj=\mathbb{E}(\tau^+_j);~~\tmj=\mathbb{E}(\tau^-_j);~~\pi_i=\Pn(\lambda)(\{i\}).}
Applying Lemmas~2.1 and 2.2 of \cite{BX01} with birth rate $\lambda$, death rate $\beta_i=i$, $A:=[k,\infty)$ and $\pi(\cdot)=\sum_{l\in \cdot}\pi_l$, we have 
\ben{\label{sol2}f(i)=\tmi\pi(A\cap[0,i-1])-\tpi1\pi(A\cap[i,\infty)),~~i\ge 1} 
and for $j\in\Z_+$,\ben{\label{notation1}\tpj=\frac{F(j)}{\lambda\pi_j},~~\tmj=\frac{\overline{F}(j)}{j\pi_j},}
where, {as in Theorem~\ref{thm1},}
 \ben{\label{notation2}F(j)=\sum_{i=0}^j\pi_i;~~\overline{F}(j)=\sum_{i=j}^{\infty}\pi_i.}
One can easily simplify \Ref{sol2} to get
 \ben{\label{sol3}f(i)=\begin{cases}
-\tpi1\pi(A)&\mbox{for }i\le k,\\
-\tmi F(k-1)&\mbox{for }i>k,
\end{cases}}
which, together with \Ref{notation1} and the balance equations
\begin{equation}
\lambda \pi_i=(i+1)\pi_{i+1},\mbox{ for all }i\in\Z_+,\label{balancedeq}
\end{equation}
implies
\ben{\label{diff11}\Delta f(i)=\begin{cases}
-\pi(A)\left(\frac{F(i)}{\lambda\pi_i}-\frac{F(i-1)}{\lambda\pi_{i-1}}\right)~~&\mbox{for }i\le k-1,\\
-(1-\pi(A))\left(\frac{\overline{F}(i+1)}{\lambda\pi_i}-\frac{\overline{F}(i)}{\lambda\pi_{i-1}}\right)~~&\mbox{for }i\ge k.
\end{cases}}
It follows from \cite[Lemma~2.4]{BX01} that for $i\ge 1$,
\be{\frac{F(i)}{F(i-1)}\ge\frac{\lambda}{i}\ge\frac{\overline{F}(i+1)}{\overline{F}(i)},} which,
together with \Ref{balancedeq}, ensures 
 \bean{\label{diff12}\Delta f(i)\le 0~~~~&\text{for~}i\le k-1,\\
\label{diff13}\Delta f(i)\ge 0~~~~&\text{for~}i\ge k,
}hence, $f(k)\le f(i)\le 0$ and combining \Ref{notation1}, \Ref{notation2} and \Ref{sol3} gives $\|f\|=\vert f(k)\vert=\frac{F(k-1)}{k\pi_k}\pi(A),$ as claimed in (i).

Apropos of (ii), because of \Ref{diff12} and \Ref{diff13}, it remains to show that $\Delta f$ is decreasing in the two ranges. To this end,
we will mainly rely on the properties of the solution \Ref{sol1}. Let $T$ be an exponential random variable with mean $1$ and independent of 
birth-death process $Z_{i-1}$, then $Z_i$ can be represented as 
\be{Z_i(t)=Z_{i-1}(t)+\bone_{\{T>t\}},} 
hence we obtain from \Ref{sol1} and the strong Markov property in the second last equality that
\bea{f(i)&=-\int_0^{\infty}\mathbb{E}\left[\bone_{\{Z_{i-1}(t)+\bone_{\{T>t\}}\ge k\}}-\bone_{\{Z_{i-1}(t)\ge k\}}\right]dt\non\\
&=-\mathbb{E}\int_0^{\infty}e^{-t}\bone_{\{Z_{i-1}(t)=k-1\}}dt\non\\
&=-\mathbb{E}\left\{\int_{\tau_{i-1,k-1}}^{\infty}e^{-t}\bone_{\{Z_{i-1}(t)=k-1\}}dt\right\}\non\\
&=-\mathbb{E}\left\{e^{-\tau_{i-1,k-1}}\right\}\mathbb{E}\int_0^{\infty}e^{-t}\bone_{\{Z_{k-1}(t)=k-1\}}dt\non\\
%&=-\mathbb{E} e^{-\tau_{i-1,k-1}}\mathbb{E}\left\{\int_0^{\infty}e^{-t}\bone_{\{Z_{k-1}(t)=k-1\}}dt\right\}\non\\
&=\mathbb{E} e^{-\tau_{i-1,k-1}}f(k).}
This enables us to give another representation of \Ref{diff11} as
 \ben{\Delta f(i)=f(k)(\mathbb{E} e^{-\tau_{i,k-1}}-\mathbb{E} e^{-\tau_{i-1,k-1}}),\label{conj}}
 and so
 \be{\Delta^2f(i)=f(k)(\mathbb{E} e^{-\tau_{i+1,k-1}}-2\mathbb{E} e^{-\tau_{i,k-1}}+\mathbb{E} e^{-\tau_{i-1,k-1}}).}
For %$i-1\ge k-1$, i.e. 
$i\ge k$, using the strong Markov property again in the equalities below, we have
\bea{&\mathbb{E} (e^{-\tau_{i+1,k-1}}-2e^{-\tau_{i,k-1}}+e^{-\tau_{i-1,k-1}})\non\\
&=\mathbb{E} e^{-\tau_{i-1,k-1}}(\mathbb{E} e^{-\tau_{i+1,i-1}}-2\mathbb{E} e^{-\tau_{i,i-1}}+1)\non\\
&=\mathbb{E} e^{-\tau_{i-1,k-1}}(\mathbb{E} e^{-\tau_{i+1,i}}\mathbb{E} e^{-\tau_{i,i-1}}-2\mathbb{E} e^{-\tau_{i,i-1}}+1)\non\\
&\ge \mathbb{E} e^{-\tau_{i-1,k-1}}(\mathbb{E} e^{-\tau_{i,i-1}}-1)^2\ge0,}
where the inequality follows from
 \bea{&\tau_{i,i-1}=\inf\{t:Z_{i}(t)=i-1\}\\
 &=\inf\{t:Z_{i}(t)+\bone_{\{T>t\}}=i-1+\bone_{\{T>t\}}\}\\
&\ge\inf\{t:Z_{i+1}(t)=i\}=\tau_{i+1,i}.}
Similarly, for $i\le k-2$, $\tau_{i-1,i}$ is stochastically less than or equal to $\tau_{i,i+1}$, so
\bea{&\mathbb{E} (e^{-\tau_{i+1,k-1}}-2e^{-\tau_{i,k-1}}+e^{-\tau_{i-1,k-1}})\non\\
&=\mathbb{E} e^{-\tau_{i+1,k-1}}(\mathbb{E} e^{-\tau_{i-1,i+1}}-2\mathbb{E} e^{-\tau_{i,i+1}}+1)\non\\
&\ge\mathbb{E} e^{-\tau_{i+1,k-1}}(\mathbb{E} e^{-\tau_{i,i+1}}-1)^2\ge 0.
} Hence, $\Delta^2 f(i)\le 0$ for $i\ge k$ and $i\le k-2$, which concludes the proof of (ii).

In terms of (iii), 
we use (ii) to obtain   
\bea{\|\Delta f\|_{k-}&=\vert\Delta f(k-1)\vert=f(k-1)-f(k)\non\\
&=\pi(A)\frac{1}{\lambda}\left(\frac{F(k-1)}{\pi_{k-1}}-\frac{F(k-2)}{\pi_{k-2}}\right)\non\\
&={\pi(A)\frac{F(k-1)}{k\pi_k}\left(1-\frac{F(k-2)}{F(k-1)}\cdot\frac{\lambda}{k-1}\right)}.}
Likewise, 
\bea{\|\Delta f\|_{k+}&=\vert\Delta f(k)\vert=f(k+1)-f(k)\non\\
&=\frac{F(k-1)}{\lambda\pi_{k-1}}\pi(A)-\frac{\overline{F}(k+1)}{\lambda\pi_k}F(k-1)\non\\
&={\pi(A)\frac{F(k-1)}{k\pi_k}\left(1-\frac{\overline{F}(k+1)}{\overline{F}(k)}\cdot\frac{k}{\lambda}\right)}.}
Since (iv) is clearly an immediate consequence of (iii), \Ref{Steinconstantnew1} and \Ref{Steinconstantnew2}, the proof of Lemma~\ref{lma2} is complete. \qed
\vskip10pt
\noindent{\bf Proof of Theorem \ref{thm1}}~~As in the proof of Lemma~\ref{lma2}, we set $A=[k,\infty)$ and $h=\bone_A$, 
then
\be{\mathbb{P}(W-a\ge k)-\mathbb{P}(Y\ge k)=\mathbb{E} h(W-a)-\Pn(\lambda)\{h\}.} Define \bea{
e_1:=&\mathbb{E}(h(W-a)-\Pn(\lambda)\{h\})\bone_{\{W-a<0\}}-\lambda f(0)\mathbb{P}(W-a=-1),\\
e_2:=&\mathbb{E}(\lambda f(W-a+1)-(W-a)f(W-a)),
}then it follows from \Ref{steinid} that \ben{\label{trans1}\mathbb{P}(W-a\ge k)-\mathbb{P}(Y\ge k)=e_1+e_2.}

For the estimate of $e_1$, from $f(0)=f(1)$, we know that $\lambda f(0)=-\Pn(\lambda)\{h\}$, thus \be{e_1=-\mathbb{P}(W-a<-1)\Pn(\lambda)\{h\},}
which gives
\ben{\label{error1}\vert e_1\vert=\pi(A)\mathbb{P}(W-a<-1).} 

For the estimate of $e_2$, denoting $\tf(j):=f(j-a)$, we have
\ben{e_2=\mathbb{E}\left\{\lambda\Delta\tf(W)-(W-\mu)\tf(W)\right\}.\label{e2-1}
}

Using Lemma~\ref{lma2} (ii), we have $\Delta^2\tf(m)$ is negative for all $m$ except $m=a+k-1$, which implies  $-\sum_{m\ne k-1}\Delta^2 f(m)\le \Delta^2 f(k-1)=\|\Delta^2 f\|$ and
$$\mean\left[\left.\Delta^2\tf(W_i'+l)\right|\cf_i\right]\le \Delta^2 f(k-1) \mathbb{P}\left[\left.W_i'=k-1+a-l\right|\cf_i\right]\le \|\Delta^2 f\|\theta_i $$
and
$$\mean\left[\left.\Delta^2\tf(W_i'+l)\right|\cf_i\right]\ge \sum_{m\ne k-1}\Delta^2 f(m) \mathbb{P}\left[\left.W_i'=m+a-l\right|\cf_i\right]\ge -\theta_i \|\Delta^2 f\|,
$$
hence
\begin{equation}
\left|\mean\left[\left.\Delta^2\tf(W_i'+l)\right|\cf_i\right]\right|\le \|\Delta^2 f\|\theta_i. \label{conditionalbound}
\end{equation}

By taking \be{\theta:=\lambda-\sigma^2,}
we have from \Ref{e2-1} that
\bean{e_2&=\theta\mathbb{E}\Delta\tf(W)+\mathbb{E}\left\{\sigma^2\Delta\tf(W)-(W-\mu)\tf(W)\right\}\non\\
&=\theta\mathbb{E}\Delta\tf(W)+\mathbb{E}\left\{\sigma^2\Delta\tf(W)-\sum_{i\in\ci}(X_i-\mu_i)\tf(W)\right\}\non\\
%&=\theta\mathbb{E}\tf(W)+\mathbb{E}\left\{\lambda\Delta\tf(W)-\sum_{i\in\ci}(X_i-\mu_i)\left(\tf(W)-\tf(W_i)\right)\right\}\non\\
&=\theta\mathbb{E}\Delta\tf(W)+\sigma^2\mathbb{E}\Delta\tf(W)-\sum_{i\in\ci}\mathbb{E}\left\{(X_i-\mu_i)\left(\tf(W)-\tf(W_i)\right)\right\}\non\\
&=\theta\mathbb{E}\Delta\tf(W)+\sigma^2\mathbb{E}\Delta\tf(W)-\sum_{i\in\ci}\mathbb{E}\left\{(X_i-\mu_i)\left(\sum_{j=0}^{Z_i-1}\Delta\tf(W_i+j)\right)\right\}\non\\
&=\theta\mathbb{E}\Delta\tf(W)+\sigma^2\mathbb{E}\Delta\tf(W)-\sum_{i\in\ci}\mathbb{E}\left[(X_i-\mu_i)Z_i\right]\mathbb{E}\Delta\tf(W_i')\non\\
&\ \ \ \ \ \ \ \ -\sum_{i\in\ci}\mathbb{E}\left\{(X_i-\mu_i)\sum_{j=0}^{Z_i-1}\left[\Delta\tf(W_i+j)-\Delta\tf(W_i')\right]\right\}\non\\
&=\theta\mathbb{E}\Delta\tf(W)+\sum_{i\in\ci}\mathbb{E}\left[(X_i-\mu_i)Z_i\right]\mathbb{E}\left[\Delta\tf(W)-\Delta\tf(W_i')\right]\non\\
&\ \ \ \ \ \ \ \ -\sum_{i\in\ci}\mathbb{E}\left\{(X_i-\mu_i)\sum_{j=0}^{Z_i-1}\left[\Delta\tf(W_i+j)-\Delta\tf(W_i')\right]\right\}\non\\
&=\theta\mathbb{E}\Delta\tf(W)+\sum_{i\in\ci}\mathbb{E}\left[(X_i-\mu_i)Z_i\right]\mathbb{E}\left[\sum_{j=0}^{Z_i'-1}\Delta^2\tf(W_i'+j)\right]\non\\
&\ \ \ \ \ \ \ \ -\sum_{i\in\ci}\mathbb{E}\left\{(X_i-\mu_i)\sum_{j=0}^{Z_i-1}\sum_{l=0}^{Z_i'-Z_i+j-1}\Delta^2\tf(W_i'+l)\right\}\non\\
&=\theta\mathbb{E}\Delta\tf(W)+\sum_{i\in\ci}\mathbb{E}\left[(X_i-\mu_i)Z_i\right]\mathbb{E}\left[\sum_{j=0}^{Z_i'-1}\mathbb{E}\left(\left.\Delta^2\tf(W_i'+j)\right|\cf_i\right)\right]\non\\
&\ \ \ \ \ \ \ \ -\sum_{i\in\ci}\mathbb{E}\left\{(X_i-\mu_i)\sum_{j=0}^{Z_i-1}\sum_{l=0}^{Z_i'-Z_i+j-1}\mathbb{E}\left(\left.\Delta^2\tf(W_i'+l)\right|\cf_i\right)\right\},\label{xiathm1-1}}
where the third last equality is because $\sum_{i\in\ci}\mathbb{E}[(X_i-\mu_i)Z_i]=\sigma^2$ and $(X_i,Z_i)$ is independent of $W_i'$, and the last equality is due to the assumption 
that $\{X_j:\ j\in B_i\}$ is $\cf_i$ measurable. Using \Ref{conditionalbound} in \Ref{xiathm1-1},
we obtain
\bean{
&|e_2|\le \|\Delta f\||\theta|\non\\
&+\|\Delta^2 f\|\sum_{i\in\ci}\theta_i\left\{{|\mathbb{E}(X_i-\mu_i)Z_i|}\mathbb{E} (Z_i')+\mathbb{E}\left[|X_i-\mu_i|Z_i(Z_i'-Z_i/2-1/2)\right]\right\}.\label{xiathm1-2}}
Now, combining Lemma~\ref{lma2} (iii), (iv), \Ref{trans1}, \Ref{error1} 
and \Ref{xiathm1-2} gives \Ref{esti1}. \qed

\vskip10pt
\noindent{\bf Proof of Corollary~\ref{cor1}} Under the setting of the local dependence, the claim follows from Theorem~\ref{thm1} by taking $Z_i=Z_i'=X_i$. \qed

\vskip10pt
\noindent{\bf Proof of Theorem \ref{thm1.2}}~~Recall the Stein representation \Ref{trans1} and the estimate \Ref{error1}, it remains to tackle \Ref{e2-1}. However,
\bea{e_2&=\mathbb{E}\left(\lambda \tf(W+1)-\mu \tf(W^s)+a\tf(W)\right)\non\\
&=\mu\mathbb{E}(\tf(W+1)-\tf(W^s))+(\lambda-\mu)\mathbb{E}\Delta \tf(W),}
thus\bean{\label{error2}
\vert e_2\vert&\le\|\Delta f\|(\mu\mathbb{E}\vert W+1-W^s\vert+\vert\lambda-\mu\vert)\non\\
&\le\Pn(\lambda)\{h\}\left[\bc_1(\lambda,k)(\mu\mathbb{E}\vert W+1-W^s\vert+\vert\lambda-\mu\vert)\right].}
Hence, combining \Ref{trans1}, \Ref{error1} and \Ref{error2} completes the proof.\qed

\vskip10pt
\noindent{\bf Proof of Theorem \ref{thm2}}~~Again, we make use of the Stein representation \Ref{trans1} and the estimate \Ref{error1} so that it suffices to deal with \Ref{e2-1}. To this end, we have
\bea{e_2&=\mathbb{E}\left(\lambda \Delta\tf(W)-(W-\mu) \tf(W)\right)\non\\
&=\mathbb{E}\left(\lambda \Delta\tf(W)-\sigma^2 \Delta\tf(W^\star)\right)\non\\
&=\mathbb{E}\left((\lambda-\sigma^2) \Delta\tf(W)+\sigma^2 (\Delta\tf(W)-\Delta\tf(W^\star))\right).}
However, with $R=W^\star-W$,
\bea{&\mathbb{E}\left[\Delta\tf(W)-\Delta\tf(W^\star)\right]\non\\
&=-\mathbb{E}\left\{\sum_{j=0}^{R-1}\mathbb{E}\left(\Delta^2\tf(W+j)\right)\bone_{R>0}-\sum_{j=1}^{-R}\mathbb{E}\left(\Delta^2\tf(W-j)\right)\bone_{R<0}\right\}\non\\
&=-\mathbb{E}\left\{\sum_{j=0}^{R-1}\mathbb{E}\left(\left.\Delta^2\tf(W+j)\right|R\right)\bone_{R>0}-\sum_{j=1}^{-R}\mathbb{E}\left(\left.\Delta^2\tf(W-j)\right|R\right)\bone_{R<0}\right\},}
and a similar argument for \Ref{conditionalbound} ensures
$$\left|\mathbb{E}\left(\left.\Delta^2\tf(W+j)\right|R\right)\right|\le \|\Delta^2f\|\theta_R,$$
hence 
\bean{\label{error4}
\vert e_2\vert&\le|\lambda-\sigma^2|\|\Delta f\|+\sigma^2\|\Delta^2f\|\mathbb{E}[|R|\theta_R].}
The claim follows from combining \Ref{trans1}, \Ref{error1} and \Ref{error4} and using Lemma~\ref{lma2} (iii), (iv).\qed

\vskip10pt
\noindent{\bf Proof of Proposition~\ref{secondpropPnBi}}~~The first inequality of \Ref{Poissonbinomial1} is a direct consequence of \cite{Hoeffding56}. For the second inequality, let $h=\bone_{[k,\infty)}$ and $f$ be the solution of the Stein identity \Ref{steinid} with $\lambda=\mu$, set $W_i=W-X_i$, $Y\sim\Pn(\mu)$, the following argument is standard (see \cite[p.~6]{BHJ}) and we repeat it for the ease of reading:
\begin{eqnarray}
&&\mathbb{P}(W\ge k)-\mathbb{P}(Y\ge k)\nonumber\\
&&=\mathbb{E}\{\mu f(W+1)-Wf(W)\}\nonumber\\
&&=\mu \mathbb{E} f(W+1)-\sum_{i=1}^n\mathbb{E} \{X_if(W)\}\nonumber\\
&&=\mu \mathbb{E} f(W+1)-\sum_{i=1}^np_i\mathbb{E} \{f(W_i+1)\}\nonumber\\
&&=\sum_{i=1}^np_i^2\mathbb{E}\Delta f(W_i+1).\label{secondpropPnBi01}
\end{eqnarray}
For any non-negative integer valued random variable $U$ such that the following expectations exist, the summation by parts gives
$$
\mathbb{E} g(U+1)=\sum_{j=1}^\infty \Delta g(j)\mathbb{P}(U\ge j)+g(1).$$
On the other hand, \cite[Proposition~2.1]{BCC95} ensures
$$\frac{\mathbb{P}(W_i\ge j)}{\mathbb{P}(Y\ge j)}\le \frac{\mathbb{P}(W\ge j)}{\mathbb{P}(Y\ge j)}\le \sup_{r\ge 0}\frac{\mathbb{P}(W=r)}{\mathbb{P}(Y=r)}\le M,$$
so using Lemma~\ref{lma2} (ii), we have
\begin{eqnarray}
\mathbb{E}\Delta f(W_i+1)&=&\sum_{j=1}^\infty \Delta^2f(j)\mathbb{P}(W_i\ge j)+\Delta f(1)\non\\
&\ge&M\sum_{j\ge 1, j\ne k-1} \Delta^2f(j)\mathbb{P}(Y\ge j)+\Delta f(1)\nonumber\\
&=&M\left\{\sum_{j=1}^\infty \Delta^2f(j)\mathbb{P}(Y\ge j)+\Delta f(1)\right\}\nonumber\\
&&-M\Delta^2f(k-1)\mathbb{P}(Y\ge k-1)+(1-M)\Delta f(1)\nonumber\\
&=&M\mathbb{E}\Delta f(Y+1)-M\Delta^2f(k-1)\mathbb{P}(Y\ge k-1)+(1-M)\Delta f(1)\nonumber\\
&>&M\mathbb{E}\Delta f(Y+1)-M\Delta^2f(k-1)\mathbb{P}(Y\ge k-1).\label{secondpropPnBi03}
\end{eqnarray}
However, by \Ref{sol1}, since $\Pn(\mu)$ is the stationary distribution of $Z_i$, $Z_{Y}(t)\sim \Pn(\mu)$, leading to 
\begin{eqnarray}
&&\mathbb{E}\Delta f(Y+1)\nonumber\\
&=&-\int_0^\infty \mathbb{E}[h(Z_{Y+2}(t))-2h(Z_{Y+1}(t))+h(Z_{Y}(t))]dt\nonumber\\
&=&-\int_0^\infty \mathbb{E}[h(Y+\bone_{\{T_1>t\}}+\bone_{\{T_2>t\}})-h(Y+\bone_{\{T_1>t\}})-h(Y+\bone_{\{T_2>t\}})+h(Y)]dt\nonumber\\
&=&-\int_0^\infty e^{-2t}\mathbb{E}[\Delta^2h(Y)]dt=-\frac12 (\pi_{k-2}-\pi_{k-1}),\label{secondpropPnBi04}
\end{eqnarray}
{where $T_1,T_2$ are i.i.d. $\exp(1)$ random variables independent of $Y$.}
Combining \Ref{secondpropPnBi01}, \Ref{secondpropPnBi03} and \Ref{secondpropPnBi04},
we have
{\begin{eqnarray}
&&\frac{\mathbb{P}(W\ge k)}{\mathbb{P}(Y\ge k)}-1\nonumber\\
&&> -\frac12M\mu_2\frac{\pi_{k-2}-\pi_{k-1}}{\mathbb{P}(Y\ge k)}-M\mu_2\Delta^2f(k-1)\frac{\mathbb{P}(Y\ge k-1)}{\mathbb{P}(Y\ge k)}.\label{secondpropPnBi05}
\end{eqnarray}}
For the first term of \Ref{secondpropPnBi05}, using \cite[Proposition~A.2.1 (ii)]{BHJ}, we obtain
\bean{
\frac{\pi_{k-2}-\pi_{k-1}}{\mathbb{P}(Y\ge k)}&=\frac k {\mu}\cdot \frac{\pi_k}{\mathbb{P}(Y\ge k)}\cdot\frac{k-1-\mu}{\mu}\non\\
&\le \frac {4(k-\mu)} {\mu}\cdot\frac{k-1-\mu}{\mu}\le 4x^2/\mu.\label{secondpropPnBi08}
}
For the second term of \Ref{secondpropPnBi05}, we use the crude estimate of {$\Delta^2 f(k-1)\le 2 \|\Delta f\|\le 2(1-e^{-\mu})/\mu$} (see \cite[Lemma~1.1.1]{BHJ} or Remark~\ref{remarknaive}), so applying \cite[Proposition~A.2.1 (ii)]{BHJ} again,
{\bean{&\Delta^2 f(k-1)\frac{\mathbb{P}(Y\ge k-1)}{\mathbb{P}(Y\ge k)}\non\\
&\le \frac{2(1-e^{-\mu})}\mu\left(1+\frac{\pi_k}{\mathbb{P}(Y\ge k)}\cdot \frac k\mu\right)\non\\
&\le  \frac{2(1-e^{-\mu})}\mu\left(1+\frac{4(k-\mu)}{\mu}\right)\le \frac2{\mu}\left(1+4x\sqrt{\frac{1-e^{-\mu}}\mu}\right).\label{secondpropPnBi09}
}}
The bound \Ref{Poissonbinomial1} follows by collecting \Ref{secondpropPnBi05}, \Ref{secondpropPnBi08} and \Ref{secondpropPnBi09}.
\qed

\vskip10pt
\noindent\textbf{Acknowledgements} We thank the anonymous referees for suggesting the ``naive bound'' in Remark~\ref{remarknaive} and 
comments leading to the improved version of the paper. We also thank Serguei Novak for email discussions about the quality of the bounds presented in the paper versus the
``naive bound''. 

%%%%%%%%%%%%%%%%%%%%%%%%%%%%%
%%%%%%% References   %%%%%%%%
%%%%%%%%%%%%%%%%%%%%%%%%%%%%%

\def\ac{{Academic Press}~}
\def\aap{{Adv. Appl. Prob.}~}
\def\ap{{Ann. Probab.}~}
\def\anap{{Ann. Appl. Probab.}~}
\def\eljp{{\it Electron.\ J.~Probab.\/}~} 
\def\jap{{J. Appl. Probab.}~}
\def\jws{{John Wiley~$\&$ Sons}~}
\def\ny{{New York}~}
\def\ptrf{{Probab. Theory Related Fields}~}
\def\sp{{Springer}~}
\def\spa{{Stochastic Process. Appl.}~}
\def\sv{{Springer-Verlag}~}
\def\tpa{{Theory Probab. Appl.}~}
\def\zw{{Z. Wahrsch. Verw. Gebiete}~}

\end{document}